%\magnification=1200

\def\e{\epsilon}  \def\g{\gamma}         
\def\s{\sigma}  \def\S{\Sigma}

\def\cl{{\rm {cl}}}                 \def\HD{{\rm{HD}}}

  \def\a{\alpha}        \def\La{\Lambda}
\def\z{\zeta}             \def\k{\kappa}
\def\e{\varepsilon}      \def\b{\beta}
\def\f{\varphi}            \def\d{\delta}
\def\la{\lambda}   \def\th{\vartheta}

\def\C{I\!\!\!\!C}     \def\D{I\!\!D}    \def\Z{Z\!\!\!Z}
\def\R{I\!\!R}          \def\T{{\cal T}}   
\def\Crit{{\rm {Crit}}}      
      \def\N{{\cal N}}

\def\and{{\rm and}}   \def\ov{\overline}   
 \def\dist{{\rm {dist}}}     
\def\h{{\rm h}}                   \def\supp{{\rm {supp}}}
        \def\diam{{\rm {diam}}}
\def\Const{{\rm {Const}}}   \def\deg{{\rm {deg}}}
   \def\ess{{{\rm {ess}}}}
\def\diaml{\diam_{\Re\log}}   \def\Tel{{\rm {Tel}}}

{\bf Accessability of  typical points for invariant measures of positive Lyapunov exponents for iterations of holomorphic maps}

\

\centerline {by F. Przytycki\footnote{*}{ supported by Polish KBN 
Grants 210469101 "Iteracje i Fraktale" and 210909101 
"...Uklady Dynamiczne".} }

\

%\centerline {(version of  April 20, 1993)}

\

{\bf Abstract.}  {\it We prove that  if A is the basin of 
immediate attraction to a periodic attracting or parabolic 
point for a rational map f on the Riemann sphere, 
if $A$ is completely invariant (i.e. $f^{-1}(A)=A$), and if 
$\mu$ is 
an arbitrary 
 $f$-invariant  measure with positive Lyapunov 
exponents
on 
$\partial A$, then $\mu$-almost every point 
$q\in\partial A$ is 
accessible along a curve from $A$. In fact we prove the accessability of every "good" $q$ i.e. such $q$ for which 
"small 
neighbourhoods arrive at large scale" under iteration of 
$f$.

This generalizes Douady-Eremenko-Levin-Petersen 
theorem on the accessability of periodic sources. 

We prove a general "tree" version of this theorem. This 
allows to deduce  that on the limit set of a geometric coding tree (in particular on the whole 
Julia set), if diameters of the edges converge to 0 
uniformly with the number of generation converging to 
$\infty$, every $f$-invariant probability ergodic  
measure with positive Lyapunov exponent is the image 
through coding 
with the help of the tree, of an invariant measure on 
 the full one-sided shift space.

The assumption that $f$ is holomorphic on 
$A$, or on the domain $U$ of the tree, can be  relaxed and 
one does not need to assume $f$ extends beyond $A$ or 
$U$.

Finally we prove that in the case $f$ is polynomial-like on 
a neighbourhood of $\ov\C\setminus A$ every 
"good" $q\in\partial A$ is accessible along an external 
ray.}

\

{\bf Introduction.} 

\

Let $f:\ov\C\to \ov\C$ be a rational map of the Riemann sphere $\ov\C$. Let $J(f)$ denote its Julia set. 
We say a periodic point $p$ of period $m$ is attracting (a sink) 
if $|(f^m)'(p)|<1$,
repelling (a source) if $|(f^m)'(p)|>1$
 and parabolic if $(f^m)'(p)$ is a root of 
unity. We say that $A=A_p$ is the immediate basin of 
attraction to a sink or a parabolic point $p$ if $A$ is a component of \ $\ov\C\setminus 
J(f)$ such that $f^{nm}|_A\to p$ as $n\to\infty$ and 
$p\in A_p$ in the case $p$ is attracting, $p\in\partial A$
in the case $p$ is parabolic.

\

We call $q\in\partial A$ {\it good} \  if there exist 
real numbers $r>0, 
\k >0, \d : 0<\d<r$ and an integer $\Delta>0$ such that for every $n$ large enough 
$$
\sharp\{\ \hbox{good times} \ \} /n \ge\k  \eqno (0.0)
$$ 
We call here 
${\ov n}: 0 \le {\ov n} \le n$ a {\it good time} if 
for each $0\le l\le {\ov n}-\Delta$ the component 
$B_{{\ov n}, l}$ of $ f^{-({\ov n}-l)}(B(f^{{\ov n}}(q),r)$ containing 
$f^l(q)$ satisfies: 

$$B_{{\ov n},l}\subset B(f^l(q),r-\d) \eqno (0.1)$$

In the definition of {\it good} $q$ we assume also that 
$$
\lim_{{\ov n}\to\infty} \diam (B_{{\ov n},0}) \to 0   \eqno (0.2)
$$
lim taken over {\it good} ${\ov n}$'s.

Finally in the definition of {\it good} $q$  we assume 
about each {\it good} $\ov n$ that 

$$f^{-{\ov n}}(A)\cap B_{{\ov n},0} \subset A .   \eqno (0.3)$$

\

\

We shall prove the following 

{\bf Theorem A.} Every good $q\in\partial A$ is 
accessible from $A$, i.e. 
 there exists a continuous 
curve  $\g:[0,1]\to \ov\C$ such that $\g([0,1))\subset 
A$  and $\g(1)=q$.

\

Theorem A generalizes Douady-Eremenko-Levin-Petersen theorem on the accessability of periodic 
sources. Remark that in the case of  periodic sources one 
obtains curves 
along which periodic $q $ is accessible, of finite lengths, 
see Section 1. 
Condition (0.1) holds in the case $q$ is a periodic source 
 for all ${\ov n}$'s. Condition (0.3) 
is  true if $A$ is the basin of attraction to 
$\infty$ for $f$ a polynomial, and more generally if 
$A$ is completely invariant, i.e $f^{-1}(A)=A$. 

Condition (0.3) in the case of a source is equivalent 
to Petersen's condition [Pe]. 

\

Under the assumption of the complete invariance of $A$ 
\ $\mu$-almost every point for $\mu$ an invariant 
probability measure with positive Lyapunov exponents is {\it good} hence accessible, cf. Corollary 0.2.

\

In fact  we shall introduce in Section 2 a weaker 
definition of {\it good} $q$ and prove Theorem A 
with that weaker definition. In that weaker definition parabolic periodic points in $\partial A$ are good. 
The traces of {\it telescopes } built there can sit in an 
arbitrary interpetal, so one obtains the accessability in 
each interpetal. One obtains in particular Theorem 18.9 
in [Mi1].

Remark that the above conditions of being {\it good} 
are already quite weak.  In particular 
we do not exclude critical points in 
$B_{{\ov n},l}$.  

For example every point in $\partial A$ 
is {\it good} if $A$ is the basin of attraction to $\infty$ for a polynomial $z\mapsto z^2+c$ which is non-renormalizable, $c$ outside the "cardioid".  This is 
Yoccoz-Branner-Hubbard theory, see [Mi2]. (In this case 
however theorem A is worthless because
one proves directly the local connectedness of $\partial A$.)

\

Remark that complete invariance of $A$,  a basin of attraction to a sink, 
does 
not imply that $f$ on a neighbourhood of 
$\ov\C\setminus A$ is polynomial-like. (Polynomial-like maps were first defined and studied in [DH].)
In [P4]
an example of degree 3, of the form $z\to z^2+c+{b\over z-a}$, 
with a completely invariant
basin of attraction to $\infty$, not simply-connected, with
only 2 critical points in the basin, is described.

 \

We prove in the paper a theorem more general 
than  Theorem A, namely a theorem on the accessability along  
branches  of a {\it geometric coding tree}.  We  
recall now basic definitions from  [P1, P2, PUZ, PS]. 

Let $U$ be an open connected subset of the Riemann sphere $\ov\C$.
Consider any holomorphic mapping $f:U\to \ov\C$ such that $f(U)\supset U$ and $f:U\to f(U)$ is a proper map. 
Denote $\Crit(f)=\{z:f'(z)=0\}$. This is called the set of 
critical  points for $f$. Suppose that $\Crit(f)$ is finite.
 Consider any $z\in f(U)$. Let $z^1,z^2,...,z^d$ be all
 the $f$-preimages of $z$ in $U$ where $d=\deg f\ge 2$. 
(Pay attention that we consider here, unlike in the other papers, only the full tree i.e. not only some preimages but 
all preimages of $z$ in $U$.)

Consider smooth curves $\g^j:[0,1]\to f(U)$, \ $j=1,...,d$,   joining $z$ with $z^j$ respectively (i.e. $\g^j(0)=z, \g^j(1)=z^j$), such that there are no critical values for iterations of $f$ in $\bigcup _{j=1}^d \g^j $, i.e. 
$\g^j\cap f^n(\Crit(f))=\emptyset$ for every $j$ and $n>0$.  We allow 
self-intersections of each $\g^j$.

Let $\S^d:=\{1,...,d\}^{\Z^+}$ denote the one-sided shift space and $\s$ the shift to the left, i.e. $\s((\a_n))=(\a_{n+1})$. 
We consider the standard metric on $\S^d$
$$
\rho((\a_n),(\b_n))=\exp -k((\a_n),(\b_n))
$$
where $k((\a_n),(\b_n))$ is the least integer for which 
$\a_k\not=\b_k$.

For every sequence $\a=(\a_n)_{n=0}^\infty \in \S^d$ we define $\g_0(\a):=\g^{\a_0}$. Suppose that for some $n\ge 0$, for every $ 0\le m\le n$, and all $\a\in\S^d$, the curves $\g_m(\a)$ are already defined.
Suppose that for $1\le m\le n$ we have $f\circ \g_m(\a)=\g_{m-1}(\s(\a))$, and $\g_m(\a)(0)=\g_{m-1}(\a)(1)$.

 Define the curves $\g_{n+1}(\a) $ so that the previous equalities hold by taking respective $f$-preimages of curves $\g_n$. For every  $\a\in\S^d$ and $n\ge 0$ denote $z_n(\a):=\g_n(\a)(1)$. 

For every $n\ge 0$ denote by $\S_n=\S^d_n$ the space of all 
sequences of elements of $\{1,...,d\}$ of length $n+1$. 
Let $\pi_n$ denote the projection $\pi_n:\S^d\to 
\S_n$ defined by $\pi_n(\a)=(\a_0,...,\a_n)$. As 
$z_n(\a)$ and $\g_n(\a)$ depends only on 
$(\a_0,...,\a_n)$, we can consider $z_n$ and $\g_n$ as functions on $\S_n$.

The graph $\T=\T(z,\g^1,...,\g^d)$ with the vertices $z$ and $z_n(\a)$ and 
edges $\g_n(\a)$ is called a {\it geometric coding tree} with the root at 
$z$. For every $\a\in\S^d$ the subgraph composed of $z,z_n(\a)$ and 
$\g_n(\a)$ for all $n\ge 0$ is called a  {\it geometric branch} and denoted 
by $b(\a)$. The branch $b(\a)$ is called {\it convergent} if the sequence 
$\g_n(\a)$ is convergent to a point in $\cl U$. We define the {\it coding map} 
$z_\infty :{\cal D}(z_\infty)\to \cl U$ by $z_\infty(\a):=\lim_{n\to\infty}
z_n(\a)$ on the domain ${\cal D}={\cal D}(z_\infty)$ of all such $\a$'s for 
which $b(\a)$ is convergent.

\

In Sections 1-3, for any curve (maybe with self-intersections) $\g:I\to\ov\C$
where $I$ is a closed interval in $\R$, we call 
$\g$ restricted to  $J$ a subinterval (maybe degenerated to a point) of $I$ 
a {\it part} of $\g$.
Consider $\g$ on $J_1\subset [0,1]$ and $\g'$ on $J_2\subset [0,1]$ either
both $\g$ and $\g'$ being 
parts of one $\g_n(\a)$, $J_1\cap J_2=\emptyset, \ J_1$ between 0 and 
$J_2$ , or $\g$ a part of $\g_{n_1}(\a)$ and $\g'$ a part of $\g_{n_2}$ where
$n_1< n_2$. Let $\Gamma:[0,n_2-n_1+1]\to\ov\C$ be
the concatenation of $\g_{n_1},\g_{n_1+1},...,\g_{n_2}$. We call the  
restriction of $\Gamma$ to the convex hull 
of $J_1\subset [0,1]$ and $J_2\subset [n_2-n_1, n_2-n_1+1]$ (we identified here
$[0,1]$ with $[n_2-n_1, n_2-n_1+1]$) {\it a part of $b(\a)$ between $\g$ and 
$\g'$ }.

\

For every continuous map $F:X\to X$ of a compact 
space $X$ denote by  $M(F)$ the set of all probability 
$F$-invariant measures on $X$.
In the case $X$ is a compact subset of the Riemann 
sphere $\ov\C$ and $F$ extends holomorphically to a  
neighbourhood of $X$ and $\mu\in M(F)$ we can 
consider for $\mu$-a.e. $x$ Lyapunov characteristic 
exponent
$$
\chi(F,x)=\lim_{n\to\infty}{1\over n}\log |(F^n)'(x)|
.$$
If $\mu$ is ergodic then for $\mu$-a.e. $x$ 
$$
\chi(F,x)=\chi_\mu(F)=\int\log(F')d\mu
.$$

In this paper where we shall discuss properties of 
$\mu$-a.e. point, it is enough to consider only ergodic  
measures, because by Rochlin Decomposition Theorem 
every $\mu\in M(F)$ can be decomposed into ergodic ones.

Denote
$$
M^{\chi+} _e(F)=\{\mu\in M(F) : \mu \ \hbox{ergodic} \  
\chi_\mu(F)>0\}
$$
$$
M^{\h+} _e(F)=\{\mu\in M(F) : \mu \ \hbox{ergodic} \  
\h_\mu(F)>0\}
$$
where $\h_\mu(F)$ denotes measure-theoretic entropy.

   From Ruelle Theorem it follows that
$\h_\mu(F) \le 2\chi_\mu(F)$ see [R], so $M^{\h+}_e(F) 
\subset M^{\chi+}_e(F) $.

\

The basic theorem concerning convergence of 
geometric coding trees  is the following:

\

{\bf  Convergence Theorem.} 1. Every branch except branches in a set of Hausdorff dimension  0 in the  metric $\rho$ on $\S^d$, is 
convergent. (i.e $\HD(\S^d\setminus{\cal D} )=0$).
In particular for every $\nu\in  M^{\h+}(\s) $ we have 
$\nu(\S^d\setminus{\cal D} )=0$, so the measure 
$(z_\infty)_*(\nu)$ makes sense.

2. For every $z\in\cl U$, \ $\HD(z_\infty^{-1}(\{z\}))=0$. 
Hence for every $\nu\in M(\s)$ we have for the 
entropies:  
$\h_{\nu_\f}(\s)=\h_{ (z_\infty)_*(\nu_\f)}(\ov f)> 0$,
(if we assume that there exists $\ov f$ a 
continuous extension of $f$ to $\cl U$).

\

The proof of this Theorem can be found in 
[P1] and [P2] under some  assumptions on a slow convergence of $f^n(\Crit(f)$ to $\g^j$ for $n\to\infty$)  
and in [PS] in full generality ( even with $f^n(\Crit(f))\cap\g^j\not=\emptyset$  allowed).

\

Let $\hat\La$ denote the  set of all limit points of $f^{-
n}(z), n\to\infty$. Analogously to the case $q\in \partial A$ we say that $q\in\hat\La $ is {\it 
good} if $f$ extends holomorphically to a neighbourhood 
of $\{f^n(q), n=0,1,...\}$ ( we use the same symbol $f$ to 
denote the  extension) and  conditions (0.0'), (0.1'),  (0.2') 
and (0.3') hold. These conditions are defined similarly 
 to (0.0)-(0.3), with $A$ replaced by $U$ and $\partial  A$
replaced by $\hat\La$.

Again pay attention that we shall give a precise weaker definition 
of $q$ {\it good} in Section 2. and prove Theorem B with that weaker definition. That definition will not demand 
$f$ extending beyond $U$.

\

{\bf Theorem B.}  Let $f$ be  a holomorphic mapping 
$f:U\to \ov\C$ and $\T$ be a geometric coding tree in 
$U$ as above. Suppose
$$
\diam (\g_n(\a)) \to 0  \ \hbox{as} \  n\to 
\infty      \eqno(0.4)
$$
uniformly with respect to $\a\in\S^d$.

Then every good $q\in\hat\La$ is a limit point of a 
branch $b(\a)$. 

\

Using a lemma belonging to  Pesin Theory 
(see Section 2) we prove that $\mu$-a.e.$q$ below is good and easily obtain the following 

{\bf Corollary 0.1.} Let $f$ be  a holomorphic mapping 
$f:U\to \ov\C$ and $\T$ be a geometric coding tree in 
$U$ such that the condition (0.4) holds. 
  If $\mu$ is a probability measure on 
$\hat\La$ and  the map $f$ extends holomorphically 
from $U$ to a neighbourhood of $\supp\mu$ so that 
$\mu\in M_e^{\chi+}(f)$,
then for $\mu$-almost every $q\in\hat\La$ satisfying (0.3')  
there exists $\a\in\S^d$ such that $b(\a)$ converges to 
$q$. In particular $\mu$ is 
a $(z_\infty)_*$-image of a measure $m\in M(\s)$ on $\S^d$.

\

Remark that Corollary 0.1 concerns in particular every 
$\mu$ with  $\h_\mu(f)>0$.  Assuming 
that $f$ extends holomorphically to a neighbourhood of 
$\hat\La$ and refering also to Convergence 
Theorem we see
that $(z_\infty)_*$ maps 
 $M^{\h+}_e(\s)$ onto $M^{\h+}_e(f|_{\hat\La})$ preserving entropy.   

The question whether this correspondence is {\it onto} 
is stated in [P3].  Thus Corollary 0.1 answers this 
question in positive under additional assumptions (0.3') 
and (0.4).

We do not know whether this correspondence is finite-to-one except measures supported by orbits of periodic sources  for which the answer is positive, see Proposition 1.2.  

\

Two special cases are of  particular interest. The first 
one  corresponds to Theorem A:

\

{\bf Corollary 0.2.} Let $f:\ov\C\to\ov\C$ be a rational 
mapping and $A$ be a completely invariant basin 
of attraction to a sink or a parabolic point. Then for every 
$\mu\in M_e^{\chi+}(f|_{\partial A})$ \ $\mu$-a.e. 
$q\in\partial A$ is accessible  from $A$.

\

{\bf Corollary 0.3.}  Let $f:\ov\C\to\ov\C$ be a rational 
mapping, $\deg f=d$, and $\T=\T(z,\g^1,...,\g^d) $ be a 
geometric coding tree.  Assume (0.4). Let $\mu\in 
M_e^{\chi +}(f)$. Then for $\mu$-a.e. $q$   
there exists $\a\in\S^d$ such that $b(\a)$ converges to 
$q$.

\

\

  In Theorem A and Corollary 0.2 in the case 
$f$ is a polynomial (or a polynomial-like map)  
and $A$ is the basin of attraction to 
$\infty$, the accessability of a point along a curve often implies automathically the accessability along an external ray. In the case $A$ is simply-connected this follows from Lindel\"of's Theorem. 
External 
rays are  defined as images under standard 
Riemann map of rays $t\z, \z\in \partial\D, 1<t<\infty$.

In the case $A$ is not simply-connected one should first define  external rays in the absence of 
Riemann  map.  This is done in [GM] and [LevS] in the 
case of $f$ a polynomial and in [LevP] in the polynomial-
like situation. We recall these definitions in Section 3. 

We prove in Section 3 the following

\

{\bf Theorem C.} Let $W_1\subset W$ be open, connected, simply-connected domains in $\ov\C$ 
such that $\cl W_1 \subset W$ and $f:W_1\to W$ be a 
polynomial-like map. denote $K=\bigcap_{n\ge 0}f^{-
n}(W)$. Then every good $q\in \partial K$ is accessible along an external ray in $W\setminus K$.

\

An alternative way to prove the accessability along an external ray is to use somehow, as in the simply-connected case, Lindel\"of's Theorem. This is performed 
in [LevP]. It is proved there that if $q$ is accessible along a curve in $W\setminus K$ and $q$ belongs to a periodic or preperiodic component $K(q)$ of $K$ then 
it is accessible along an external ray.

Pay attention also that for any $q\in\partial K$ if $K(q)$ 
is one point then $q$ is accessible along an external ray.
This is easy, see [GM, Appendix] and [LevP].

\

\

{\bf Remark 0.4. (Proof of Theorem A from B and Corollary 0.2 from 0.1).}  We do not know how to get rid of the assumption (0.4) in Theorem B and Corollary 1.  
In Theorem A and Corollary 2 this condition is 
guaranteed  automathically. 
More precisely to deduce Theorem A from B and Corollary 2 from 1 we consider an arbitrary tree
$\T=\T(z,\g^1,...,\g^d) $ in $A$, where $d=\deg(f|_A)$,
 so that $\g^j\cap\bigcup_{n>0}f^n(\Crit(f))=\emptyset$ 
and $p\notin\bigcup_{j=1,...,d}\g^j$. Only critical points 
in $A$ account here. Forward orbits of these critical 
points  converge to $p$ hence the following condition holds:
$$
\Bigl(\bigcup_{j=1,...,d}\g^j \Bigr)
\cap \cl \Bigl(\bigcup_{n>0} 
f^n(\Crit(f))\Bigr) = \emptyset                        \eqno (0.5)
$$

Hence we can take open discs $U^j\supset\g^j$ 
such that
$$
\bigcup_{j=1,...,d}U^j \cap \cl \Bigl(\bigcup_{n>0} 
f^n(\Crit(f))\Bigr) = \emptyset 
$$
and consider univalent branches $F_n(\a)$ of $f^{-n}$ 
mapping respective  $\g^j$ to $\g_n(\a)$. \ 
$\{F_n(\a)\}_{\a,n}$ is 
a normal family of maps.  If it had a non-constant limit 
function $G$ then we would find an open domain $V$ 
such that $F_{n_t}(V)\subset U$ as $n_t\to\infty$.  If we assumed $p\notin U^j$ we arrive at a contradiction.
This proves (0.4). Finally by the complete invariance 
of $A$ we have $\hat\La=\partial A$.

\

In Corollary 0.3 to find $\T$ such that (0.4) holds it is 
enough to assume that the forward limit set of 
$f^n(\Crit (f))$ does not dissect $\ov\C$, because then 
we find $\T$ so that (0.5) holds.

We believe however that in Proof of 
Corollary 3 we can omit (0.4), 
or maybe often find a tree such that (0.3) holds.

\

{\bf Remark 0.5.} Observe that there are examples where 
(0.4) does not hold. Take for example $z$ in a Siegel disc or $z$ being just a sink.   Even if $J(f)=\ov\C$ one should be careful: for M. Herman's examples 
$z\mapsto \la z{z-a\over 1-\ov a z}/{z-b \over 1-\ov b 
z}, \ |\la|=1, a\not=0\not= b, a\approx b$,see [H1], the unit circe is invariant and for a branch in it (0.4) 
fails.  These examples are related with the notion of neutral sets, see [GPS].

\

{\bf Remark 0.6.} The assumption $f$ is holomorphic on $U$ (or $A$) can be replaced by the assumption $f$ is 
just a continuous map, a branched cover over  $f(U)\supset 
U$. 

However without the holomorphy of $f$ we do not 
know how the assumption (0.4) could be verified. 

\

{\bf Remark 0.7.} The fact that in, say, Theorem A we do 
not need to assume that $f$ extends holomorphically 
beyond the basin $A$ suggests that maybe the assumption (0.3) is substantial 
and without it the accessability in Theorem A  is not true.  We have in mind 
here an analogous situation of a Siegel disc with 
the boundary not simply-connected, where the map is only smooth beyond it, see [H2]. Accessability of  periodic 
sources in the boundary of $A$ in the absence of the 
assumption (0.3) is a famous open problem and we think 
that if the answer is positive one should substantially use 
in a proof the holomorphy of $f$ outside $A$.

\

The paper is organised as follows: in Section 1 we prove 
theorem B for $q$ a periodic source, in Section 2 we deal 
with the general case. The case of sources was known in the polynomial-like and parabolic $p$ situations
[D], [EL], [Pe]. The general case contains the case of sources but it  is  more tricky (though not more 
complicated)  so we decided to separate the case of sources to make the paper more understandable.  
Section 3 is devoted to Theorem C.

\

\

\

{\bf Section 1.  Accessability of periodic sources.}

\

{\bf Theorem D.}  Let $f:U\to\ov\C$ be a holomorphic 
map and ${\cal T}(z,\g^1,...\g^d)$ be a geometric 
coding tree in $U$, \  $d=\deg f|_U$. 
Assume (0.4).
Next assume  that $f$ extends holomorphically  to a neighbourhood of a family 
of points $q_0,...,q_{n-1} \in \hat\La$ so that this family 
is a periodic repelling orbit for this extension (the 
extension is also denoted by $f$). 

Assume  finally that there exists $V$  a neighbourhood of 
$q$ on which $f^n$ is linearizable and if $F$ is its inverse 
on $V$ such that $F(q)=q$ then 
$$
F(V\cap U) \subset U   \eqno (1.1)
$$

Then there exists a periodic $\a\in\S^d$ such that 
$b(\a)$ is convergent to $q$. Moreover the convergence 
is exponential, in particular the curve being the body of 
$b(\a)$ is of finite length. 

\

{\bf Proof of Theorem D.} As usually we can suppose 
that $q$ is a fixed point by passing to the iterate $f^n$ if 
$n>1$.  

Assume that $q\not= z$. We shall deal with the case 
$q=z$  later.

\

Let  $h$ denote the linearizing map  i.e. a map conjugating $f$ 
on a neighbourhood of  
$\cl V$ to $z\to \la z$ with $\la=f'(q)$, mapping $q$ to 
$0\in\C$. 

Replace if necessary the set $V$ by a smaller 
neighbourhood of $q$ so that $z\notin V$ and 
 $\partial V =h^{-1}\exp\{\Re\xi=a\}$
for a constant $a\in \R$.

\

For every set $K\subset \cl V \setminus 
\{q\}$ consider its diameter in the radial direction 
(with origin at $q$) in the logarithmic scale, namely 
the diameter  of the projection of the set $\log h(K)$ to the 
real axis.  This  will be denoted by $\diam_{\Re\log}(K)$.

\

For every $m\ge 0$  write 
$$
R_m:= h^{-1}\exp(\{\z\in\C: 
a-(m+1)\log|\la|  < \Re\z < a-m\log|\la| \})
$$
and
$$
V_m:=h^{-1}\exp(\{\z\in\C: 
 \Re\z < a-m\log|\la| \}).
$$

\

Observe the following important property of $\g_n(w)$'s, 
$n\ge 0, w\in\S^d$ : 

For every $\e>0$ there exists $N(\e)$ such that 
if a component $\g$ of $\g_n(w) \cap R_m$ satisfies  
$$
\diam_{\Re\log} (\g) > \e\log|\la|
 \ \hbox{and} \ z_n(w)\in V_m               \eqno (1.2)
$$
then 
$$
0 < n-m < N(\e)                                              \eqno (1.3)
$$

Indeed,   by (1.2)
for every $t=0,1,...,m$  we have 
  $f^t(z_n(w))\in V_{m-t}$ so $f^t(z_n(w))\not=z$.
 Hence $n>m$.
On the other hand 
 we have 
$$
\e\le\diam_{\Re\log}(\g)=\diam_{\Re\log}(f^{m}
(\g)\le \Const  \ \diam (f^{m}(\g)
$$
So from (0.3) and from the estimate
$\diam f^m(\g_n(w))=\diam \g_{n-m} 
(\s^m(w))  \ge \e$, we deduce that $n-m$
 is bounded by a constant depending 
only on $\e$. This proves (1.3).

\

Fix topological discs $U^1,...,U^d$ being 
neighbourhoods of $\g^1,...,\g^d$ respectively 
such that $\bigcup_{i=1}^{N(\e)}  f^i(\Crit (f)) \cap U^j
=\emptyset$ for every $j=1,...,d$. 

(There is a minor 
inaccuracy here because this concerns the case the curves $\g^j$ are embedded. If they have self-intersections we should cover them by families of small 
discs and later lift them by branches of $f^{-t}$ one by one along the curves.)

For every $\g$ being a part of  $\g_n(w)$ satisfying (1.2) we can 
consider 
$$W_1=F_{n-(m-1)}(\s^{m-1}(w))(U^j)$$ 
which is a 
neighbourhood of $f^{m-1}(\g_n(w))$. We used here 
the notation $F_t(v)$ for the branch of $f^{-t}$ mapping 
$\g^j$ to $\g_t(v), v\in \S^d$. Here $j=v_t$.

Next consider the component $W_2$ of $W_1\cap V$ 
containing $f^{m-1}(\g)$. Using Koebe's Bounded 
Distortion 
Theorem we can find a disc 
$$
W(\g)=B(x, \Const \e \la^{-m})        \eqno     (1.4)
$$ 
in $F^{m-1}(W_2)$ with $x\in\g$ such that $f^n$ 
maps  $W(\g)$ univalently into $U^j$. We take $\Const$ 
such that 
$$
\diaml W(\g)  <  {1 \over 2}\log |\la|.       \eqno  (1.5)
$$
(Remark that this part is easier if (0.5) is assumed. Then 
we just consider $U^j$'s disjoint with 
$\cl\bigcup_{n=1}^\infty \Crit(f)$.)

\

By the definition of $\hat\La$ 
there exist $n_0\ge 0$ and  $\a \in \S^d$ 
such that $\g_{n_0}(\a) \cap V \not= \emptyset$ . By (1.1) there 
exist $\b_1,\b_2,... $ each in $\{1,...,d\}$ such that 
for each $k \ge 0$ we have 
$$
F^{k}(b(\a))=b(\b_k,\b_{k-1},...,\b_1,\a).
$$
More precisely we consider an arbitrary component 
$\hat\g$ of 
$\g_{n_0}(\a) \cap V$ and extend $F^k$ from it holomorphically along $b(\a)$.

Denote for abbreviation $\b_k,\b_{k-1},...,\b_1,\a$ by $k]\a$.

Denote also $F^k(\hat\g)$ by $\hat\g_{k]}$ and the part 
of $\g_{n_0+k}(k]\a)$ between $\hat\g_{k]}$ and
$z_{n_0+k-1}(k]\a)$ by $\g_{k]}$.

\

For each $k\ge 0$ denote by ${\cal N}_k$ 
the set of 
all pairs of integers $(t,m)$ such that $t: 0\le t \le k+n_0, 
0<m<k$ and $\g_t(k]\a)$ 
satisfies (1.2) for a curve $\g$ being a part of $\g_t(k]\a)$ 
\  or a part of $\g\subset\g_{k]}$ if $t=k+n_0$ 
and for the integer $m$ 
and additionally 
$$  \{ \ \hbox{the part of} \ b(\a) \ \hbox{between} \ \g \ \hbox{and} \ 
\hat\g_{k]} \}\subset V_m .               \eqno (1.6)                     
$$ 
         
We write in this case $W(\g)=W_{k,t,m}$ and 
$\g=\g_{k,t,m}$.  Figure 1 illustrates our definitions:

\

\

\

\

\

\

\

\

\

\

\

\

\centerline {Figure 1.}

\

We have now two possibilities:

1. For every  
$k_2>k_1\ge 0, \  0<m_1<k_1, 0<m_2<k_2$ and 
$0\le T \le k_2+n_0$ such that $(T,m_1)\in {\cal 
N}_{k_1}, \  (T,m_2)\in {\cal N}_{k_2}$,
supposed the equality of the $T$-th entries 
$(k_1]\a)_T = (k_2]\a)_T$, we have 

$$W_{k_1,T,m_1} \cap W_{k_2,T,m_2} =\emptyset.$$

(The equality of the $T$-th entries means  that
$f^T(W_{k_1,T,m_1}), f^T(W_{k_2,T,m_2})$ are in  the 
same $U^j$.)

\

2.  The case 1. does not hold,  what 
implies obviously the existence of $T$ and the other 
integers as 
above such that  $\pi_T(k_1]\a)=\pi_T(k_2]\a)$, 
(i.e. the blocks of $k_1]\a$ and $k_2]\a$ from 0
to $T$ are the same).

\

Later we shall prove that the case 1. leads to a contradiction. Now we shall prove that the case 2. 
allows to find a periodic branch convergent to $q$ 
what proves our Theorem.

\

Denote $K=k_2 - k_1$. Repeat that we have 
$$
\pi_T(\s^K(k_2]\a)) = \pi_T(k_1]\a) = \pi_T(k_2]\a).
$$
Denote $k_2]\a$ by $\th$. We get by the above:
$$
f^K(z_{T+K}(\th))=z_T(\th).
$$
or writing this with the help of $F$ which is the inverse 
of $f$ on $V$ so that $F(q)=q$ we have 
$F^K(z_T(\th))=z_{T+K}(\th).$
We know also that $\g:=\bigcup_{t=T+1}^{T+K} 
\g_t(\th)$ being a curve joining 
$z_T(\th)$ with $z_{T+K}(\th))$   is contained in $V$
 (even in $V_{m(k_2,t)}$) by (1.4).

Hence the curve $\Gamma:=\bigcup_{n\ge 
0}F^{nK}(\g)$ is the body of the part starting from 
the $T$-th vertex of the periodic branch
$(\th_0,...,\th_{K-1},\th_0,...,\th_{K-1},\th_0,...)$.

\

To finish Proof of Theorem D we should now 
eliminate the disjointness case 1.  We shall 
just prove there is not enough room for that.

\

Denote for every $k\ge 0$
$$
A^+_k:=\{m:  0< m < k, \ \hbox{there 
exists }\  t  \  \hbox{such that} \  (t,m)\in {\cal N}_k \}
$$
Let $A^-_k:=\{1,...,k-1\} \setminus A^+_k$.

\

 As $\g_{k+n_0}(k]\a) $ intersects $V_k$ (at 
$\hat\g_{k]}$),  
each $0< m \le k-1$ is fully intersected by the 
curve built from the curves $\g_t(k]\a), t=0,...,k+n_0-1$ 
and $\g_{k]}$.

Hence
$$
\sharp A^-_k \log|\la|    
 \le \sum_{m\in A_k^-}   \Bigl(\sum_{0\le s \le n_0+k}
\diam_{\Re\log} (\g_s(k]\a)\cap R_m) \Bigr)
\le 2(k+n_0+1) \e\log|\la|.
$$

The coefficient 2 takes into account the possibility that 
one $\g_s(k]\a)$ intersects $R_m$ and $R_{m+1}$, 
where $m,m+1\in A^-_k$ (it cannot intersect more than 
two $R_m$'s because $\diaml (\g_s(k]\a)\cap R_m) <\e$).

\

Hence
$$
\sharp A^-_k \le 2(k+n_0+1)\e.
$$
So
$$
\sharp A^+_k \ge k-2(k+n_0+1)\e-1  \ge k(1-3\e) \eqno 
(1.7)
$$
for $k$ large enough.

\

Fix from now on $\e=1/4$. Fix an arbitrary large $k_0$. 
Let ${\cal N}^+=\bigcup_{0\le k\le k_0}(k,{\cal N}_k)$.

Observe that each point $\xi\in V$ 
belongs to at most 
$$
4dN(1/4)                                                         \eqno (1.8)
$$ 
sets $W(k,t,m)$ where $(k,(t,m))\in
{\cal N}^+$. 

Indeed 
if $W(k_1,t_1,m_1) \cap 
W(k_2,t_2,m_2)\not=\emptyset$ 
then  $|m_1 -m_2| \le 1$ by (1.5), and by (1.3) we have 
$$
|m_i-t_i|< N(1/4), \ \ i=1,2
$$
hence  
$$
|t_1 - t_2| < 2 N(1/4).
$$
(In the case $t_i=k_i+n_0$ for $i=1$ or 2 we cannot in 
fact refer to (1.3). The trouble is  with its $n-m>0$ part, because we do not know whether $z_{k_i+n_0}\in 
V_{m_i}$. But 
then directly $m_i<k_i\le t_i$.)

\

But we assumed (this is our case 1.) 
that for every $t,m$ and $j$ all the sets 
$W(k,t,m)$ with the $t$-th entry of $k]\a$ equal to $j$, variable $k$, 
are pairwise disjoint. This finishes the proof of the 
estimate (1.8).

\

The conclusion from (1.8)  and (1.4) is that because of 
the  lack of room $\sharp {\cal N}^+ <\Const 
k_0$. This contradicts (1.7) for $\e=1/4$ and $k_0$ large enough.

The disjointness case 1. is 
eliminated. Theorem D in the case $z\not=q$ is 
proved. 

\

Consider the case $z=q$. Then, unless $\g^j\equiv q$ in 
which case Theorem is trivial, the role of $z$ in the 
above proof can be 
played by arbitrary $z^j\in \g^j\setminus\{q\}$. Formally on the level 0 we have now $d^2$ curves joining each 
$z^j$ with preimages of $z^i$ in $\g_1((i,j))$. \hfill$\clubsuit$

\

{\bf Remark 1.1.} Under the assumption $z\not=q$ and moreover $q\notin\bigcup_{j=1,...,d}
\g^j $ (which is the case when we apply Theorem B to prove theorem A) observe
that there exists a constant  $M$ such  that for every $n\ge 0$ and 
$\th\in\S^d$ we have $\diaml \g_n(\th) < M$.

Indeed let $m=m_1\ge 0$ be the smallest integer such that 
$\g_n(\th)$ intersects $R_m$  and let $m_2$ be the largest one.  Suppose that $m_2-m_1>1$. 
Then by (1.3) $n<m_1+1+N(1)$ and $m_2<n$. (The role 
of $z_n(\th)$ in the proof of this part of (1.3) is played by 
$V_{m_2}\cap\g_n(\th)$.) Thus $m_2-m_1<N(1)$.

This observation allows to modify (simplify) slightly Proof of Theorem B. One
does not need (1.6) then .

\

{\bf Proposition 1.2.} Every branch $b(\a)$ 
convergent to a periodic source $q$ is periodic (i.e $\a$ 
is periodic). There is only a finite number of $\a$'s such that $b(\a)$ converges to $q$.

\

{\bf Proof.} Suppose $z\not=q$ and $b(\a)$ converges 
to $q$. We can take $V$, a neighbourhood of 
$q$, arbitrarily small. Then the constant $n_0$ will 
depend  on it. However the above proof shows that
we obtain the equality 
$$
\pi_T(k_1]\a)=\pi_T(k_2]\a)
$$
for $k_1-k_2$ bounded by a constant independent 
of $n_0$. $z\not= q$ implies that $T\to\infty$ as $V$ shrinks to $q$. So there exists a finite block of symbols 
$\b$ such that $\a=\b\b\b...\b\a'$ \ ($\a'$ infinite) with arbitrarily many $b$'s. So $\a$ is periodic. This 
consideration  gives also a bound for the period of $\a$ 
hence it proves finitness of the set  of  $\a$'s with $b(\a)$ convergent to $q$. \hfill $\clubsuit$

\

Remark that with some additional effort we could obtain 
an estimate for the number of branches convergent to 
$q$. In the case $q$ is in the boundary of a basin of 
attarction to a sink this estimate should give so called 
Pommerenke-Levin-Yoccoz inequality (see for example 
[Pe]).

\

\

{\bf Section 2. Theorem B and Corollary 0.1.}

\

Given $f:U\to\ov\C$ a holomorphic map and ${\cal T}=\T (z,\g^1,...,\g^d)$ a geometric coding tree in $U$  as in Introduction we shall give a definition of $q\in\hat\La$ {\it 
good} more general then in Introduction.

Let us start with some preliminary definitions:

\

{\bf Definition 2.1.} $D\subset U$ is called $n_0$-{\it significant} if  there exists $\a\in \S^d$
and $0\le n\le n_0$  such that $\g_n(\a)\cap D\not=
\emptyset$.

\

{\bf Definition 2.2.} For every $\d,\k>0$ and integer 
$k>0$  a pair of sequences 
$(D_t)_{t=0,1,...,k}$ and $(D_{t,t-1})_{t=1,...,k}$ is called a {\it telescope} or 
a $(\d,\k,k)$-{\it telescope} if each $D_t$ is an open connected subset of $U$, there exists a strictly  
increasing  sequence of 
integers $0=n_0, n_1,...,n_k$ such that
each $D_{t,t-1}$ is a nonempty component of $f^{-(n_t-
n_{t-1})}(D_t)$ contained in $D_{t-1}$ (of course $f^{n_t-
n_{t-1}}$ can have critical points in $D_{t,t-1}$),

$$t/n_t > \k \ \ \hbox{for each} \ \ t ,    \eqno (2.0)$$
and else 
$$
\dist (\partial_U^\ess D_{t,t-1}, \partial_U D_{t-1})>\d. \eqno (2.1)
$$
 Here the subscript $U$ means the boundary in $U$ and 
the essential boundary $\partial_U^\ess D_{t,t-1}$ 
is defined as $\partial_U D_{t,t-1}\setminus  
\bigcup_{n=1}^{n_t-n_{t-1}}f^{-n}(\partial U)$.

\

{\bf Definition 2.3.} A $(\d,\k,k)$-telescope is called $n_0$-{\it significant} if $D_k$ is $n_0$-significant.

\

{\bf Definition 2.4.} For any $(\d,\k,k)$-telescope we can 
choose inductively sets $D_{t,l}$, where $l=t-2, t-3,...,0$ 
by choosing $D_{t,l-1}$ as a component of 
$f^{-(n_l-n_{l-1}}(D_{t,l})$ in $D_{t-1,l-1}$. We call the 
sequence 
$$D_{k,0}\subset D_{k-1,0}\subset ... \subset D_{1,0}
\subset  D_0
$$
a {\it trace } of the telescope.

\

{\bf Definition 2.5.} We call $q\in \hat\La$ {\it good} if 
there exist $\d,\k>0$ an integer $n_0\ge 0$  and  
a sequence of $n_0$-significant 
telescopes $\Tel^k, k=1,2,...$, where $\Tel^k$ is a 
$(\d,\k,k)$-telescope,  with traces 
$D^k_{k,0}, D^k_{k-1,0},...,D^k_0$ respectively 
(to the notation of each object related to the telescope 
$\Tel^k$ we add the superscript $k$) such that

$$ D^k_{l,0} \to q  \ \ \hbox{as} \ l\to\infty \ \ \hbox{uniformly over} \ k .  \eqno (2.2)$$ 

\

{\bf Remark 2.6.} $q\in\hat\La$ good in the sense of 
Introduction (conditions (0.0')-(0.3') satisfied) is of course 
good in the above sense. Indeed we choose each 
$\Delta$'s good time and denote these times by 
$n_0,n_1,...$,  of course then $\k$ in (2.0) is $\k/\Delta$ for the old $\k$ from (0.0).

For each $k$ we define a telescope $\Tel^k$ by taking 
as $D^k_k$ an arbitrary $n_0$-significant component 
of $B(f^{n_k}(q),r)$. Such a component exists with $n_0$ 
depending only on $r$ because the set all vertices of the tree $\T$ is by definition dense in $\hat\La$.  Then 
inductively for each $0\le t<k$ we choose as $D^k_t$ 
a component of $B(f^{n_t}(q),r)\cap U$ containing 
a component $D^k_{t+1,t}$ of $f^{-(n_{t+1}-
n_t)}(D^k_{t+1})$, (such a component $D^k_{t+1,t}$ exists by (0.3'). By (0.2') an arbitrary choice of traces will 
be OK.

 Of course in the case of $U=A$ a basin of immediate attraction to a sink or a parabolic point 
one can build  telescopes with $D^k_{t,l}$ not 
containing critical points, but there is no reason for that 
to be possible in general.

\

{\bf Proof of Theorem B.}  Let $q\in \hat\La$ be a {\it 
good} point according to the definition  above. 
Fix constants $\d,\k$ and $n_0$ and a sequence of 
$\d,\k,k$ telescopes and their traces ,$k=0,1,...$ as in Definition 2.5.

We can suppose that $z\notin D^k_0$ or at least that 
each $\g^j, j=1,...,d$ has a point outside $D^k_0$. If it is 
not so then either there exists $l$ such that each $\g^j$ 
has a point outside $D^k_{l,0}$ for every $k$ in which case in the considerations below we should
consider $m\ge l$ rather than $m> 0$ or else there 
exists 
$j$ such that $\g^j\equiv q$ in which case obviously
$b(j,j,j,...)$ converges to $q$.

\

Denote $D^k_{m,0}\setminus D^k_{m+1,0}$ by 
$R^k_m$ for $m=0,1,...,k-1$ and $D^k_{k,0}$ by 
 $R^k_k$. These sets replace rings from Section 1. 

Choose for each $k$ a curve $\g_{n(k)}(\a^k)$ for 
$\a^k\in\S^d$ and $n(k)\le n_0$ intersecting $D^k_k$.
 Choose a part  $\hat\g^k$ 
of  $\g_{n(k)}(\a^k)$ in this intersection.

As in Section 1 there exists $k]\a^k=\b_0^k \b_1^k ... 
\b_{n^k_k-1}\a^k\in\S^d$  such that 
$\g_{n(k)+n^k_k}(k]\a^k)$ intersects 
$D^k_{k,0}$ and moreover it contains a part 
 $\hat\g_{k]}$ which is a lift of 
$\hat\g^k$ by $f^{n_k}$. Denote the part of
$\g_{n(k)+n^k_k}(k]\a^k)$ between $z_{n(k)+n^k_k-1}(k]\a^k)$ 
and $\hat\g_{k]}$ by $\g_{k]}$

\

Fix an integer $E>0$ to be specified later. 

Define 
${\cal N}_k$ as the set of such pairs 
$(t,m)$ that $0< m<k, 0\le t\le n^k_k+n(k)$, there exist 
integers  $E_1,E_2\ge 0, \ E_1+E_2< E$ 
 such that $\g_{t+E_2}(k]\a^k)\cap 
R^k_{m+1} \not=\emptyset$, \ 
$ \g_{t-E_1}(k]\a^k)\cap R^k_{m-1}\not=\emptyset$  and
there exists a part \ $\g(t,m)$  of  
$\g_t(k]\a^k)$  in $R^k_m$, or of  $\g_{k]}$ if $t=n^k_k+n(k)$,  
such that
 $$
 \{ \ \hbox{the part of} \ b(k]\a) \ \hbox{between} \  \g(t,m) \ \hbox{and} 
 \ \hat\g_{k]} \} \subset
  D^k_{m,0}                  \eqno (2.3)
  $$

analogously to (1.6), see Figure 2.

\

\

\

\

\

\

\

\

\centerline {Figure 2}

\

We claim that analogously to the right hand inequality of (1.3) we 
have for $(t,m)\in{\cal N}_k$
$$
t\le  n^k_{m+1} +E+ N(\d/E)          \eqno  (2.4)
$$
where $N(\e):=\sup\{n: \ \hbox{there exists} \ 
\a\in\S^d \ \hbox{such that} \ \diam(\g_n(\a)\ge \e\}
$. (The number $N(\e)$ is finite by (0.4).)

Indeed, denote the part of the curve being the 
concatenation of $\g_l(k]\a^k), l=t-E_1,...,t+E_2$ 
in $R^k_m$ joining $R^k_{m-1}$ with $R^k_{m+1}$ by 
$\Gamma$. suppose that $t-E_1 \ge n^k_m$ (otherwise the claim   
is proved). Then $f^{n^k_m}(\Gamma)$ joins a point 
$\xi\in\partial_{U}D^k_{m+1,m}$
in a curve 
$$f^{n_m^k}(\g_t(k]\a^k)=\g_{t'+n_m^k}(\s^{n_m^k}(k]\a^k)), \  \  t-E_1\le t' \le t+E_2
$$   
with 
$\partial D^k_m$.

If $\xi\notin\partial_{U}^\ess D^k_{m+1,m}$, then
$$
t' < n^k_{n_{m+1}}
$$
Otherwise there exists $n\le n^k_{m+1}$ such that $f^n(\xi)\in\partial U$. This is 
already 
outside $U$ so the trajectory of $\xi$ hits $\bigcup \g^j$ 
before the time $n^k_{m+1}$ comes.

If $\xi\in\partial_{U}^\ess D^k_{m+1,m}$ then by (2.1) 
at least one of the curves $f^{n^k_m}(\g_l(k]\a^k), 
 t-E_1\le l \le t+E_2$ has the diameter not less than 
$\d/E$. Hence 
$$
l-n_m \le N(\d/E).
$$
In both cases (2.4) is proved. 

\

 Define 
$$
A_k^+ := \{m: 0< m<n^k_k, \ \hbox{there exists} \ t \  
\hbox{such that} \ (t,m)\in {\cal N}_k\}
$$ 
and 
$$
A^-_k:= \{1,...,k-1\} \setminus A_k^+\ .
$$

As each set $R_m^k$ for $m\in A_k^-$ is crossed by a part of $b(k]\a^k)$ 
between $\g^j$ for respective $j$ and $\hat\g_{k]}$ consisting of at least $E$ edges and one edge cannot serve for more then two $R^k_m$'s 
we obtain similarly to (1.7): 
$$
E\cdot\sharp A^-_k \le 2(n^k_k +n(k)+1)
$$

Hence using (2.0) we obtain
$$
\sharp A^+_k \ge k-1-{2\over E}(n(k)+n^k_k +1) \ge 
 k(1-{3\over E\k})  \eqno (2.5)
$$

Fix from now on $E>3/\k$ and denote $\eta=1-{3\over 
E\k} >0$.

\

For every $0<M\le k$ define 
$$
A^+_k(M):=\{m\in A^+_k: m < M\}
$$

We claim that there exists $M_0>0$, not depending 
on $k$  such 
that for every $M\ge M_0, \ M\in A_k^+$ we have 
$$
\sharp A_k^+(M) \ge \eta M .    \eqno (2.6)
$$
This means that the property (2.6) true for 
$M=k$, see (2.5), extends miraculously to 
every $M\in A^+_k$ large enough. The 
proof of this claim is the same as for $A^+_k$: 

Indeed $M\in A^+_k$ implies the existence of $t$ 
such that $(t,M)\in\N_k$. By (2.3) $t\le 
n^k_{M+1}+E+N({\d\over E})$. Next we   
estimate  $\sharp A_k^+(M)$ similarly as we estimated 
$\sharp A_k^+$ with $n^k_k+n(k)+1$ replaced by 
$n^k_{M+1}+E+N({\d\over E})$. We succeed for all 
$M$ large enough.

\

Now we can conclude our Proof of Theorem B:
Let $M_n:=({1\over 2}\eta)^{-n} M_0$.
By (2.6) for every $k\ge 0$ and $n\ge 0$ there exists $m\in A^+_k$ 
such that $M_n\le m < M_{n+1}$.

For each $n=0,1,... $ there is only a finite number of 
blocks of symbols of the form $\pi_t(k]n_k)$ such that 
$(t,m)\in\N_k, \ m<M_{n+1}$.  This is so by 
(2.4). 

So there are constants $t_0\ge 0$ and 
${\cal D}_0\in\S_{t_0}$ and an infinite set
$$ \eqalign {
K_0=&\{k\ge 0: \ \hbox{there exists} \ m \ \hbox{such 
that} \ \cr 
& M_0 \le m <M_1, \ (t_0,m)\in 
\N_k ,\   
\pi_{t_0}(k]\a^k)={\cal D}_0\} \cr}
$$

In $K_0$ we find an infinite $K_1$ 
etc. by induction. For every $n>0$ we obtain
infinite $K_n\subset K_{n-1}$ and constants $t_n, {\cal 
D}_n$ 
such that 

$$\eqalign{
K_n=&\{k\in K_{n-1}: \ \hbox{there exists} \ m \ \hbox{such 
that} \ \cr & M_n \le m <M_{n+1}, \ (t_n,m)\in \N_k,\   
\pi_{t_n}(k]\a^k)={\cal D}_n\}
\cr}$$

 For $\a\in\S^d$ such that $\pi_{t_n}(\a)=\pi_{t_n}(k]\a^k)$, 
we have that $b(\a)$ converges to $q$. 

\

We assumed here that $t_n\to\infty$  as $n\to\infty$. If $\sup 
t_n=t_*<\infty$ then also ${\cal D}_n$ stabilize at ${\cal D}_*$ and 
by (2.3) $z_{t_*}({\cal D}_*)=q$. Moreover there exists a sequence of 
integers 
$j_1,j_2,...\in \{1,...,d\}$ such that $\g_t({\cal D}_*,j_1,j_2,...)\equiv q$ 
for all $t\ge t_*$ so  $b({\cal D}_*,j_1,j_2,...) $ converges to $q$.

(This is not an imaginary case. Consider a source $f(q)=q\in U$ and a tree 
$\T(q,\g^1,\g^2)$ such that  $\g_1\equiv q$ and $\g_2$ joins $q$ with 
$q'\in f^{-1}(q), q\not=q'$. Then the above proof gives $b(2,1,1,...)$ the 
branch for which  $\g_n((2,1,1,...))=\equiv q'$ for every $n\ge 1$.
\hfill$\clubsuit$

\

\

{\bf Remark 2.1.} It is curious that we did not need in the above proof neither the left hand side 
inequality (1.3): $t\ge m-\Const$ for $(t,m)\in \N_k$, 
nor the sets $W(k,t,m)$.  As mentioned already in Introduction no distortion estimates , i.e. no holomorphy 
was needed. The holomorphy of $f$ is useful only to 
verify (0.4).

\

{\bf Proof of Corollary 0.1.} This follows immediately from 
Theorem B and the following fact belonging to Pesin Theory: 

\

Let $X$ be a compact subset of $\ov\C$ and $F$ be a 
holomorphic mapping on a neighbourhood of $X$ such that $F(X)=X$.  Let $\mu\in M_e^{\chi +}(F)$. Let 
$(\tilde X, \tilde F, \tilde\mu)$ be a natural extension 
(inverse limit) of $(X,F,\mu)$. Denote by $\pi$ the projection to the 0 coordinate, $\pi:\tilde X\to X$ and by $\pi_n$ the projection to an arbitrary $n$-th coordinate. 

Then for $\tilde\mu$-a.e. $\tilde x\in\tilde X$ there exists 
$r=r(\tilde x)>0$ such that  univalent branches 
$F_n$ of $F^{-n}$ on $B(\pi(x),r)$ for $n=1,2,...$ such that 
$F_n(\pi(x))=\pi_{-n}(x))$, exist.
Moreover for an arbitrary $\l : \exp (-\chi_\mu)
  <  \lambda < 1$ (not depending on $\tilde x$) and a constant $C=C(\tilde x)>0$ 
$$
|F_n'(\pi(x))| < C\lambda^n \ \ \  \and  \ \ \ 
{|F_n'(\pi(x))\  \over  |F_n'(z)|} <C
$$
for every $z\in B(\pi(x),r), \ n>0$, 
(distances  and derivatives in the Riemann metric on 
$\ov\C$).

Moreover $r$ and $C$ are measurable functions of $\tilde x$.

\

To prove Corollary (0.1) observe that  the above fact implies the existence of numbers $r, C >0$ and a set of 
positive measure $\tilde\mu$: 
$\tilde Y\subset\tilde X$ such that the above properties 
hold for every $\tilde x\in\tilde Y$ and for these $r$ and 
$C$. Ergodicity of $\mu$ implies ergodicity of 
$\tilde\mu$. So by Birkhoff Ergodic Theorem there exists 
a set $\tilde Z\subset\tilde X$ of full measure 
$\tilde\mu$ such that for each point $\tilde x \in\tilde Z$ its 
forward orbit by $\tilde F$ hits $\tilde Y$ at the positive 
density number of times. These are {\it good times} and 
$\pi(\tilde x)$ is a {\it good} point in the sense of 
Introduction (provided they satisfy (0.3')). \hfill $\clubsuit$

\

\

{\bf Section 3. \ External rays.}

\

Let $W_1\subset W$ be open, connected, simply-connected bounded domains in the complex plane$\C$ 
such that $\cl W_1 \subset W$.  Let $f:W_1\to W$ be a holomorphic proper map "onto" $W$ of degree $d\ge 
2$. We call such a map $f$ a 
{\it polynomial-like map}. Denote $K=\bigcap_{n\ge 
0}f^{-n}(W)$. This set $K$ is called a filled in Julia set 
[DH].
We can assume that $\partial W$ is smooth. Let $M$ be an arbitrary smooth function on a neighbourhood of $\cl 
W\setminus W_1$ not having critical points, such 
that $M|_{\partial W}\equiv 0$ and $M|_{\partial 
W_1}\equiv 1$ and $M\circ f = M-1$ whereever it makes sense. Extend $M$ to $W\setminus K$ by 
$M(z)=M(f^n(z))+n$ where $n$ is such that $f^n(z)\in 
W\setminus W_1$.

\

Fix $\tau:0<\tau <\pi$ and consider curves $\g:[0,\infty)\to \cl W\setminus K$, intersecting 
lines of constant $M$ at the angle $\tau$ , (this demands 
fixing orientations), not containing critical points for $M$ with $\g(0)\in \partial W$ and converging to $K$ as the 
parameter converges to $\infty$.
 One can change the standard 
euclidean metric on $\C$ so that $\tau$ is the right angle and think about gradient lines in the new metric.  We call 
such a line a {\it smooth} $\tau$-{\it ray}. Instead of parametrizing such a curve with the gradient flow time we 
parametrize it by the values of $M$. Limits of smooth 
$\tau$-rays are called $\tau$-{\it rays}. They can pass 
through critical points of $M$. (Such a $\tau$-ray 
enters a critical point along a stable separatrix and leaves it along an unstable one, the closest clockwise or counter-clockwise. If it hits again a critical point for the 
first time it leaves it along an unstable separatrix 
on the same side from which it came 
to the previous critical point, see [GM] and [LevP] for the 
more detailed description. See Fig. 3:

\

\

\

\

\

\

\

\

\

\

\

\centerline {Figure 3.}

\

\

{\bf Proof of Theorem C.} Divide each $\tau$-ray $\g$ into pieces $\g_n, \ n\ge 1$  
each joining $f^{-n}(\partial W)$ with $f^{-n-1}(\partial 
W)$.

One  easily proves the fact corresponding to (0.4):  
$$
{\rm {length}} (\g_n)\to 0 \ \hbox{as} \ n\to\infty   \eqno (3.0)
$$
uniformly over $\tau$-rays $\g$.

The proof is the same as that of  the implication (0.5) $\Rightarrow$ 0.4) in Remark 
0.4.  We have univalent branches of $f^{-k}$ for all $k$ on neighbourhoods of $\g_n$ for external rays $\g$, neighbourhoods  not depending on $k$, for n large enough, because then critical points of $f$ in $W\setminus K$ do not interfere. There is finite number of them and  their forward trajectories escape out of $W$. 

For our $q$ find {\it significant} telescopes $\Tel^k$ as 
in Section 2, where $n_0$-significant means here 
that $D^k_k$ intersects $\g^k_{n(k)}$ for a $\tau$-ray 
$\g^k$ and $n(k)\le n_0$ a constant independent of 
$k$. This is possible by (3.0).

Denote by $\g^{k]}$ the $\tau$-ray containing a point 
of $f^{-n_k}(\g^k)$ being in $D^k_{k,0}$.

We consider, similarly to  Section 2, (2.3), 
the set
$$
{\cal N}_k=\{(t,m): \ \hbox{the same conditions as in Section 2, in 
particular} \ \g_t^{k]}\cap R_m\not=\emptyset
\}$$
Similarly we define $A_k^+$ and $A_k^+(M), M\le k$ .

The same miracle that 
$$
\sharp A^+_k(M)\ge\eta M
$$
takes place for $M\ge M_0, \  M\in A^+_k$.

To get it we prove and use the estimate $t\le m+$Const for $(t,m)\in {\cal N}_k$.

Because for $M_0\le m< M_1, (t,m)\in {\cal N}_k$ integers $t$ are 
uniformly bounded over all $k$ by say $T_0$, the parts 
$\g'^{k]}=\bigcup_{l=1}^{T_0} \g^{k]}_l$
of respective $\tau$-rays $\g^{k]}$ have a convergent subsequence and 
the limit ray (joining levels 0 and $T_0$) intersects $R_m$.

Choosing consecutive subsequences we find a limit ray converging to $q$. \hfill$\clubsuit$

\

\

{\bf References.}

\

[D] A. Douady, Informal talk at the Durham Symposium, 1988.

\

[DH] A. Douady, J. H. Hubbard,  On the dynamics of 
polynomial-like mappings, Ann. Ec. Norm. Sup. 18 (1985), 
287-343.

\

[EL] A. E. Eremenko,  G. M. Levin,  On periodic points of polynomials. Ukr. Mat. Journal 41.11 (1989), 1467-1471.

\

[GM] L. R. Goldberg, J. Milnor, Fixed points of polynomial 
maps. Part II. Fixed point portraits. Ann. scient \'Ec. Norm. Sup., $4^e$
s\'erie, 
26 (1993), 51-98.

\

[GPS] P. Grzegorczyk, F. Przytycki, W. Szlenk,
On iterations of Misiurewicz's rational maps on the Riemann sphere, Ann. Inst. Henri Poincar\'e, Physique 
Th\'eorique, 53.4 
(1990), 431-444.

\

[H1] M. Herman, Exemples de fractions rationnelles 
ayant une orbite dense sur la sph\`ere de Riemann, 
Bull. Soc. Math. France 112 (1984), 93-142.

\

[H2] M. Herman, Construction of some curious 
diffeomorphism of the Riemann sphere, J. London Math. Soc. 
34 (1986), 375-384

\

[LevP] G. Levin, F. Przytycki, External rays to periodic 
points,  manuscript.

\

[LevS] G. Levin,  M. Sodin, Polynomials with 
disconnected  Julia sets and Green maps, 
preprint 23 (1990/1991), the Hebrew University of 
Jerusalem.

\

[Mi1] J. Milnor, Dynamics in one complex variable: 
Introductory lectures, preprint SUNY at Stony Brook, IMS, 1990/5.

\

[Mi2] J. Milnor, Local connectivity of Julia sets: Expository lectures, 
preprint SUNY at Stony Brook, IMS, 1992/11.

\

[Pe]  C. L. Petersen, On the Pommerenke-Levin-Yoccoz 
inequality, preprint 

\noindent IHES/M/91/43.

\

[P1] F. Przytycki, Hausdorff dimension of harmonic 
measure on the boundary of an attractive basin for a 
holomorphic map.  Invent. Math. 80 (1985), 161-179.

\

[P2]  F. Przytycki, Riemann map and holomorphic 
dynamics. Invent. Math. 85 (1986), 439-455.

\

[P3] F. Przytycki, On invariant measures for iterations of 
holomorphic maps. In "Problems in Holomorphic 
Dynamics", preprint IMS 1992/7,  SUNY at Stony Brook.

\

[P4] F. Przytycki, Polynomials in hyperbolic components, manuscript, Stony Brook 1992.

\

[PUZ] F. Przytycki, M. Urba\'nski, A. Zdunik,  Harmonic, Gibbs and Hausdorff measures for holomorphic maps. Part 1 in Annals of Math. 130 (1989), 1-40. Part 2 in Studia Math. 97.3 (1991), 189-225. 

\

[PS] F. Przytycki, J. Skrzypczak, Convergence and pre-images of limit points for coding trees for iterations of holomorphic maps, Math. Annalen 290 (1991), 425-440.

\

[R] D. Ruelle, An inequality for the entropy of 
differentiable maps, Bol. Soc. Brasil. Mat. 9 (1978), 83-87.

\end